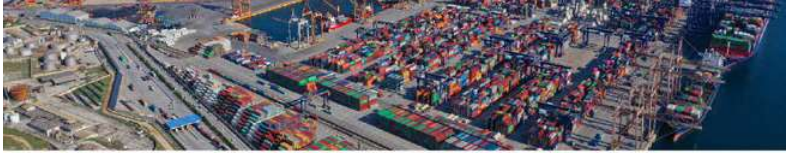

# Introducing Services and Protocols for Inter-Hub Transportation in the Physical Internet

Sahrish Jaleel Shaikh[1,2,3], Benoit Montreuil[1,2,3], Moussa Hodjat-Shamami[1,2,3], Ashish Gupta[1,2,3]

1. Physical Internet Center
2. Supply Chain and Logistics Institute
3. H.Milton Stewart School of Industrial & Systems Engineering,
Georgia Institute of Technology, Atlanta, U.S.A
Corresponding author: sahrish.shaikh@gatech.edu

## Abstract

*The Physical Internet (PI) puts high emphasis on enabling logistics to reliably perform at the speed mandated by and promised to customers, and to do so efficiently and sustainably. To do so, goods to be moved are encapsulated in modular containers and these are flowed from hub to hub in relay mode. At each hub, PI enables fast and efficient dynamic consolidation of sets of containers to be shipped together to next hubs. Each consolidated set is assigned to an appropriate vehicle so to enact the targeted transport. In this paper, we address the case where transportation service providers are available to provide vehicles and trailers of distinct dimensions on demand according to openly agreed and/or contracted terms. We describe the essence of such terms, notably relative to expected frequency distribution of transport requests, and expectations about time between request and arrival at hub. In such a context, we introduce rigorous generic protocols that can be applied at each hub so as to dynamically generate consolidation sets of modular containers and requests for on-demand transportation services, in an efficient, resilient, and sustainable way ensuring reliable pickup and delivery within the promised time windows. We demonstrate the performance of such protocols using a simulation-based experiment for a national intercity express parcel logistic network. We finally provide conclusive remarks and promising avenues for field implementation and further research.*

## 1. Introduction

The Physical Internet (PI, π) has been introduced with the specific purpose of increasing the efficiency and sustainability of logistics systems [1]. It uses the Digital Internet as a metaphor for inspiring, steering, and guiding transformations in designing, managing, and operating the supply, flow, and movement of physical objects. The Physical Internet is an open, global, hyperconnected and sustainable logistics system [2]. PI encapsulates goods in standard modular packaging, handling and transport containers. These PI containers move though distributed, multimodal transportation networks where containers from different origins are aggregated at logistic hubs depending on the next site they will travel to [3]. Open logistics facilities such as open semi-trailer transit centers and open crossdocking and consolidation hubs, as well as open warehousing, distribution, and fulfillment centers, are part of the interconnected networks, enabling a global Logistics Web [4,5]. The open logistics web is characterized by numerous logistics, handling, and buffering actors and infrastructure that are available on demand [6,7]. In this paper, we focus on hyperconnected transportation in the Physical Internet, moving modular containers from their source to through destination, through multi-segment inter-hub routes through the open mobility web.

Consider long-haul transportation as an example. Currently, long-haul drivers usually have generally a driving radius of 250 miles or more, with a maximum on the order of 11-hour trucking time per day, with variations from state to state. Going longer distances means long-haul truckers work throughout week and come home on weekends or be on the road for weeks or even months at a time. Transportation providers often have long-term contracts with shippers based on the routes their truckers will take and a fixed rate per mile. With a multi-segment network and the open mobility web such as that in the Physical Internet, the contractual terms between the shippers and transportation providers are to change drastically. On-demand trucking services have already hit the market with companies such as Convoy [8] and Uber Freight [9], offering transportation services coupled with tracking and traffic information. These are based on a pool of independent truck drivers who have registered themselves with such platforms and are able to look for shipments ready to be loaded near their vicinity. In this paper, we propose contracts between the shipper and the transportation service providers that complement such on-demand services within the Physical Internet.

Similar as in the Digital Internet, a central planning system will not be present in the Physical internet to determine how each package will move through the system. Rather, the Physical Internet can be established on a series of globally-informed and locally-applied protocols that will dynamically determine the next steps to take. The first aim of this paper is to define a set of algorithm-based transportation protocols for dynamically generating sets of modular containers to be shipped together and placing on-demand transportation service requests. These protocols are conceived for transportation in a network of open and interconnected networks where transportation service providers are available to provide vehicles of distinct dimensions on demand according to openly agreed and/or contracted terms. The second aim of the paper is to provide a simulation-based performance assessment of these protocols.

The rest of the paper is organized as follows. Section 2 focuses on the on-demand service contracts and the associated basic expectations. The design elements of the proposed protocols are described in Section 3 followed by the protocols in Section 4. The paper concludes in Section 5 by summarizing the design and motivation of the model and highlighting further research opportunities.

## 2. On-Demand Service Contracts

Freight rates are often broken down into two different categories: contractual rates and spot rates. Shippers usually sign contracts well in advance with a guaranteed volume and rates [10]. There are situations, such as inconsistent freight volumes, and seasonal or one-off shipments, in which shippers will opt for a spot rate instead [11]. However, spot rates are volatile and change with demand, and are generally significantly higher than contractual rates [10,11]. In our model, transportation service providers are available to provide vehicles and trailers of distinct dimensions, on demand according to a contracted expected frequency distribution of transport requests and order-to-arrival time of the trucks. When a package enters the system, its origin, destination, and delivery time promise are known. Based on the origin and destination of packages and modular containers, the total volume can be estimated for each route segment. The expected frequency distribution is then evaluated based on the container volume and the respective service levels expected on each segment. Important in an on-demand transportation contract is the time between the request of a vehicle and its arrival at the hub. This element is crucial as delay in the vehicle arrival can impact on-time package delivery performance. Although these times remain stochastic, the expected and latest time of arrival is included in the contract to ensure that the service level is not compromised. In the next sections, we discuss the protocol of information sharing which enables to reduce the uncertainty that revolves around expected time of arrival of vehicles.

## 3. Protocol Design Elements

In this section, we discuss the design elements involved in developing the algorithm-based protocols that dynamically generate consolidation sets of modular containers and requests for on-demand transportation services at the hubs. These elements are calculated at the beginning of a package journey and are subject to frequent update.

### 3.1. Assigned Dwell Time

The protocols consider three important location and time elements of a package; origin, destination and promised delivery time window. When a package enters the system, a route is assigned to it, which is a path from the origin to the destination, through a set of intermediary hubs. This route may be dynamically altered to deal with disruptions, yet important is maintaining a package-route duo. Based on its route, slack time is the surplus time a package has after accounting for the travel time from origin to destination and the processing time at the intermediary hubs. This excess time is distributed amongst the intermediary hubs based on the expected and relative capability and performance of the hub, and the time thus allocated is called the assigned dwell time of the package at the hub. At every intermediary node, we can calculate the latest departure time with respect to the promised delivery time and planned package route. The latest departure time of package $p$ at the first hub then is:

$$T^{ld}_{p,1} = S^{a}_{p,1} + P^{p}_{1} + D_{p,1} \quad (1)$$

For every other hub in the path,

$$T^{ld}_{p,n} = S^{a}_{p,n} + P^{p}_{n} + D_{p,n} + \delta \quad (2)$$

$$\delta = (T^{sd}_{p,n-1} - T^{ld}_{p,n-1}) \times f_n \quad (3)$$

Where $T^{ld}_{p,n}$ is the latest departure time of package $p$ from facility $n$, $S^{a}_{p,n}$ is the time of the arrival scan of package $p$ at facility $n$, $D_{p,n}$ is the assigned dwell time of the package $p$ at facility $n$, $T^{sd}_{p,n}$ is the time of the departure scan of package $p$ at facility $n$, and $f_n$ is the relative package flow of facility n. This set of formulas ensures that the latest departure time is adjusted based on the dwell time of the package in earlier hubs. If a package departed the facility earlier than its latest departure time, then it has some additional time that it can use at the next hubs. This can be also considered as an update to the assigned dwell time where we would add $\delta$ to the assigned dwell time of all remaining hubs. When calculating the assigned dwell time of a modular container, we consider all the packages in that container. The assigned dwell time of a container is the minimum of all assigned dwell times of the contents of the container. Let $\boldsymbol{P}_c$ be the set of packages in a container $c$ at hub $n$, then:

$$D_{c,n} = \min_{p \in \boldsymbol{P}_c}(D_{p,n}) \quad (4)$$

### 3.2. Package/Container Signaling

We use *package/container signaling* where the dwell time of the packages/containers is used to request the on-demand trucks. When the assigned dwell time of a package is being approached, it is flagged as an urgent one. To ensure that there is enough time to react, the contractual terms with the service provider will define an expected order-to-arrival time of the vehicle. This input considers the time it takes a truck to arrive at a facility after it has been requested. Once the urgent packages/containers are flagged, we can check the estimated time of arrival (ETA) of the packages/containers that are enroute to that facility and have the same destination as the urgent packages/containers. If the arrival time and the sorting time of incoming packages/containers is less than the remaining dwell time of current urgent packages, it means that they can be sent in the same truck. Once the truck arrives, the loading of urgent packages/containers is prioritized, and to aim to fill the truck reasonably and not below a certain threshold, the truck may be filled with all other non-urgent packages/containers destined for the same segment.

### 3.3. Information Sharing

To be well positioned to make operational-level decisions at every hub, it is crucial that pertinent data flows seamlessly across the hubs. As soon as a package is assigned a path from its origin to its destination, all the hubs that the package is planned to visit are notified and provided with the estimated time of arrival based on the path, travel times, processing times and the assigned dwell time. This information is updated when a package arrives at a hub, if it is encapsulated into a modular container, when it departs from a hub, and when any significant disruption is registered. The updated assignment of dwell time notably accounts for any instance where the package/container may have left earlier or later than planned. The hubs are also notified of the containerization status of the package, that is, whether it is travelling independently or has been consolidated with other packages into a modular container. We propose to reduce the uncertainty involved in the order-to-arrival time of the vehicles by obtaining live data on the expected arrival time of vehicles after a request has been made to the contracted vendors for each segment.

### 3.4. Maximum Latency

It may be possible to consolidate a set A of modular containers at a hub with a set B of incoming containers, yet this may require containers in set A to stay longer at the hub than their assigned dwell time. Which may lead to package delivery lateness. In order to account for this flexibility while avoiding lateness, we introduce a flexibility margin for the assigned dwell time that a container c can spend at a facility n, called Maximum Latency ($L^M_{c,n}$):

$$L^M_{c,n} = T^{ld}_{c,n} + m \left( \sum_{i \in H_c} D_{c,i} \right) \quad (5)$$

where $T^{ld}_{c,n}$ is the latest departure time of container c at facility n, based on the pre-calculated dwell time, $D_{c,i}$ is the dwell time assigned to a container $c$ at each hub $i$ element of $H_c$, the set of hubs remaining in the assigned path of container $c$. Parameter $m$ is an input between 0 and 1. Depending on the network and average number of hubs in a path, it can be altered. When $m = 0$, both the Maximum Latency and the latest departure time are same, and when $m = 1$, all dwell times are pulled to the first (current) hub. This means the package/container will have a lot of flexibility at the beginning but might be rushed at the end to meet the service level promise. As per the information sharing agreement, the hubs will receive information on the estimated time of arrival of all incoming packages; subsequently, the hubs will be able to identify any package arriving during this flexibility window.

## 4. Protocols

The protocols assume that trucking transportation service providers are available to provide vehicles and trailers of distinct dimensions on demand and that there is an established framework for sharing information. Hubs request the trucks based on the volume of containers in the hub and the incoming containers to the hub. Modular containers are flagged as urgent if they are approaching their assigned dwell time while present at a hub or if they are still en route to the hub but need to be moved as soon as they reach the hub. Hubs use these signals to request on-demand trucks, and to ensure that there is enough time to react we incorporate the contracted maximal order-to-arrival time of the vehicle into the signaling time. Therefore, if there are containers en route to a hub that require to be moved urgently, a vehicle order may be placed for them even before they arrive to the hub depending upon the order-to-arrival time of the vehicle. As the hub already has data on estimated time of arrival of the incoming containers as well as their assigned dwell time, it is able to re-evaluate the urgent containers and re-compute the queue of containers in terms of urgency. It is possible that some containers that are en route to the hub may be more urgent than the ones already present at the hub.

We hereafter discuss two protocols for creating container sets and placing vehicle requests. The first, called the Local Latency (LLT) Protocol with fixed dwell time, only considers the containers that are present at a hub and the containers that may arrive until the vehicle arrives. The second protocol, called the Maximum Latency (MLT) Protocol with flexible dwell time, proposes a flexibility margin to relax the assigned dwell time while respecting the promised delivery time window, so it also considers containers that arrive in that additional time.

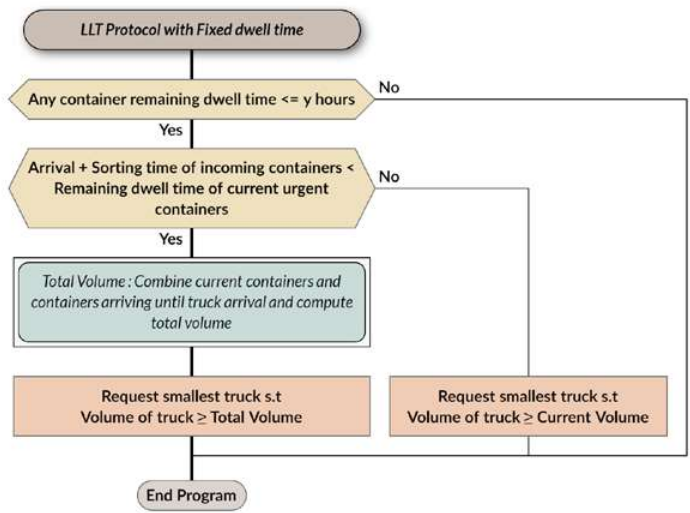

*Figure 1: Algorithm Flow Chart for Local Latency Protocol with Fixed dwell time*

### 4.1. Local Latency Protocol (LLT)

In this protocol, once the containers are flagged as urgent, the hubs verify the estimated time of arrival of incoming containers and their next destination. If there are containers that can arrive and be ready to ship within the remaining dwell time of current urgent packages, they are consolidated with the urgent packages. The facility will then place a request for the smallest number of on-demand truck(s) such that all packages can be accommodated. Once a truck arrives, it is loaded with the urgent containers, and to ensure a high fill-rate, the truck is filled with all other containers that are destined for the same segment but may not be urgent. In this paper we will use the terms trucks and trailers interchangeably. As the trailers can be of difference sizes, a truck can either be just a truck or a tractor pulling a trailer of a specific size. This algorithm ensures that the on-demand trucks are not overused, and we are able to send a lower number of larger trailers instead of higher number of small ones. The protocol is synthesized in Figure 1.

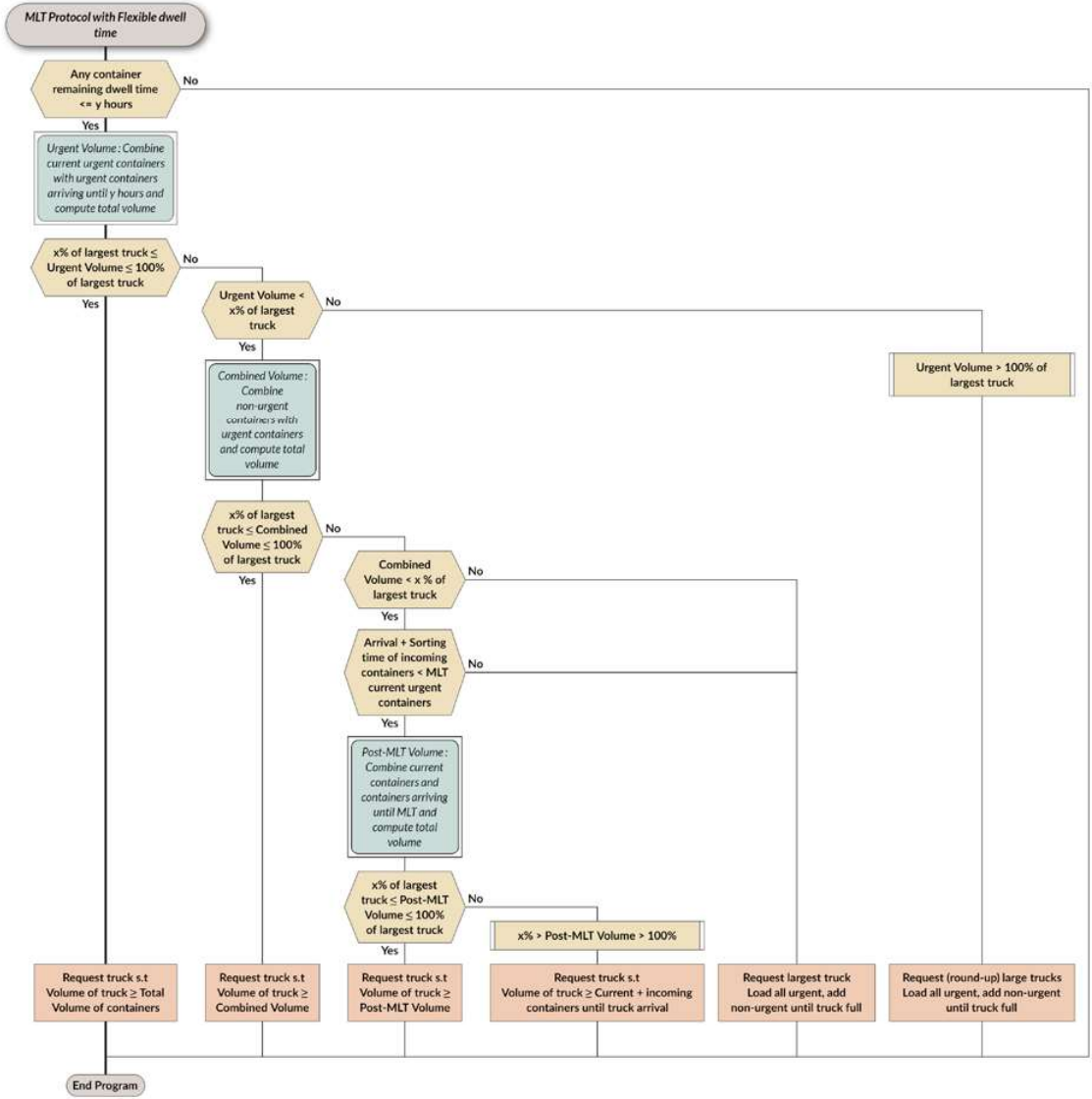

*Figure 2: Algorithm Flow Chart for MLT Protocol with Flexible dwell time*

### 4.2. Maximum Latency Protocol (MLT)

In this protocol, we introduce a flexibility margin using maximum latency. This allows the packages to wait longer at some hubs to encourage consolidation while still meeting the overall service level requirement. The hubs evaluate the volume of urgent containers currently present as well as the containers arriving within the order-to-arrival time of the vehicle. If there is not enough load to fill a large truck, the hubs check whether more containers will arrive during the Maximum Latency (MLT) of the current containers. The MLT uses a parameter $m$ between 0 and 1, calibrated to robustly ensure that the packages in the containers are to be delivered within their promised time window.

In the MLT protocol synthesized in Figure 2, input parameters $x$ and $y$ respectively refer to the minimum acceptable fill rate of a vehicle and the order-to-arrival time of a truck. In the protocol, on-demand vehicle order request is triggered as soon as the latest departure time of any package present at the hub or en route to a hub approaches the contracted expected order-to-arrival time of the vehicle. Once this protocol's algorithm is triggered, we check the volume of urgent containers available at the facility and all that may arrive and be ready to be shipped within y hours. If the volume of these urgent containers is higher than $x$, a truck is requested, and the containers are shipped. In the case where the total volume is greater than 100%, we order the minimum number of trucks required to ensure that all urgent containers move. To improve the vehicle fill rate, we fill the remaining capacity of the truck with non-urgent containers such that the ones with lower remaining dwell time amongst others are moved first. The third case is where the truck fill rate with urgent containers is below the minimum threshold. In that case, the first step is to add all non-urgent containers to evaluate if they can fill the trucks. If, after combining the two types of containers, the volume has reached the minimum threshold, we request a truck and ship the containers.

In the case where after combining non-urgent containers, the volume is more than the truck capacity, then it is already known that the volume of urgent containers was less than $x\%$, so we do not load all non-urgent containers but rather just those which can help improve the fill rate and have lower remaining dwell times than other containers. The last case is where after combining all urgent and non-urgent containers, the total volume is less than the minimum threshold. This can occur in segments with low frequency.

To check for any consolidation opportunities in this case, we use the Maximum latency concept discussed above. We check if there are any containers that may arrive if we make the current containers wait a little longer. In case, even after waiting, the expected fill rate is still under the minimum threshold, we do not wait until the MLT and request the truck immediately. However, if by including containers arriving during this flexibility window can enhance the vehicle fill rate to be above the minimum threshold, we make the current containers wait and consolidate them with the containers that arrive. The last case then is when the fill rate exceeds 100% after including the containers arriving by the MLT cutoff. If that happens, we know that there are a lot of containers arriving and those containers will suffice to request and fill a truck with satisfactory fill rate. So, we do not make the current containers wait and rather request a truck immediately to ship them.

The aim of this algorithmic protocol is to have larger trailers moving less often rather than smaller ones. If the MLT does not affect consolidation, we do not make the current containers wait to avoid any congestion and compromise on service level promises. However, if during the MLT, more consolidation opportunities are foreseen, the containers wait, and the remaining dwell times are updated for all remaining hubs for those containers. There is always a trade-off between cost and service. This protocol's algorithm can be modified by changing some parameters in order to favor cost or service.

## 5. Simulation Model
### 5.1. Framework
In order to assess the performance of the proposed protocols, we leverage a hyperconnected logistics simulator developed in Georgia Tech's Physical Internet Center, built in the AnyLogic platform according to an agent-oriented discrete-event simulation modelling approach.

Sahrish Jaleel Shaikh, Benoit Montreuil, Moussa Hodjat-Shamami, Ashish Gupta

A similar simulation program is used in [12] where the process starts with generation of customer agent either stochastically based on demand model or as deterministic historic input to the simulation. As each parcel is generated, the Inter-city agent is notified, and it makes the decision of the choice of gateway hub at origin city that the parcel needs to arrive at. The Inter-city agent decides on how to route the package as well as consolidation of packages into modular containers based on the features such as the origin, destination, direction of travel and most importantly, the committed service level or the due date of the package to the final destination [13]. The consolidation of parcels before reaching the gateway hub is done with the aim of reducing the toll on the gateway hub and facilitating the sorting/cross docking operations. Subsequently, the agent assigns the ETA, dwell times, and maximum latency for next hub and all the future hubs except the last hub. Determination of these parameters occurs on the basis of protocols discussed in the earlier sections. Intercity Router further notifies all the Gateway Hub Routers of respective Gateway Hubs in the path of incoming containers. The Gateway Hub agents manage sorting/cross docking as well as vehicle scheduling, loading/unloading operations to execute the routing decisions with the aim of minimizing operational costs within the hubs.

### 5.2. Results of Proposed Protocols

We perform the computational studies for the two algorithms using the current network infrastructure of a large China-based urban parcel delivery service provider. We focus the results on the South China region served by a network of 658 hubs, using a demand of approximately 5,000,000 packages over a time horizon of 10 days with two sizes of vehicles available. To measure efficiency, we track the total number of vehicles used, the total capacity of trucks that was on the road, the average fill rate of the trucks and the total cost incurred. To track service measures, we evaluate the service level of the packages, i.e. whether the packages reach the destination by the committed time window, and the average delay/earliness of the packages.

*Table 1: Simulation results from LLT and MLT Algorithmic Protocols*

|  | Number of vehicles | Total Trailer Capacity | Small Trailers | Large Trailers | Reduction in Trailers | Vehicle Fill Rates | Fill Rate Small Trailers | Fill Rate Large Trailers |
|---|---|---|---|---|---|---|---|---|
| **LLT Protocol** | 11,838 | 30,579,600 | 10,269 | 1,569 | 3% | 86.9% | 85.0% | 99.0% |
| **MLT Protocol** | 11,502 | 33,647,600 | 9,083 | 2,419 | | 80.4% | 79.5% | 83.7% |
|  | Total Cost ($) | Reduction in Total Cost | Service Level | Average Earliness (hours) | Sensitivity (Average earliness with tighter windows) | | | |
|  | | | | | 4 hours | 6 hours | 8 hours | 10 hours |
| **LLT Protocol** | 11,185,170 | 4% | 100.0% | 6.8 | 4.4 | 3.2 | 2.5 | 1.9 |
| **MLT Protocol** | 10,764,339 | | 100.0% | 10.0 | 7.2 | 5.8 | 3.9 | 2.1 |

Table 1 summarizes the results from both protocols. The service metrics look similar in both models as all packages reach the destination within the committed service time. We observe that although the total induced capacity from the trucks ordered increases in the MLT protocol, the total number of vehicles is reduced. This happens as MLT encourages more consolidation, so we are able to send the same packages in a lower number of larger vehicles than more frequent smaller vehicles. It is notable that both protocols perform well and all packages are delivered within the service window. That being said, we observe that the total cost is reduced by 4% in the MLT protocol. In this instance, packages move faster in the MLT protocol with average earliness of 10 hours compared to the Local Latency Protocol at 6.8 hours. We perform a sensitivity analysis on the performance of the protocols by reducing the service time windows of packages such that the promised window for delivery tightens while maintaining the feasibility of moving a package from the origin to the destination in that window. Although we make some packages wait longer in the MLT protocol, we are able to move a lot of non-urgent packages faster when consolidating. However, with very tight service levels, the earliness converges between the two protocols as the MLT protocol faces an even tighter window towards the end of the package journey once it uses the flexibility margin at the intermediate stops. We report the corresponding sensitivity results in the table above as well.

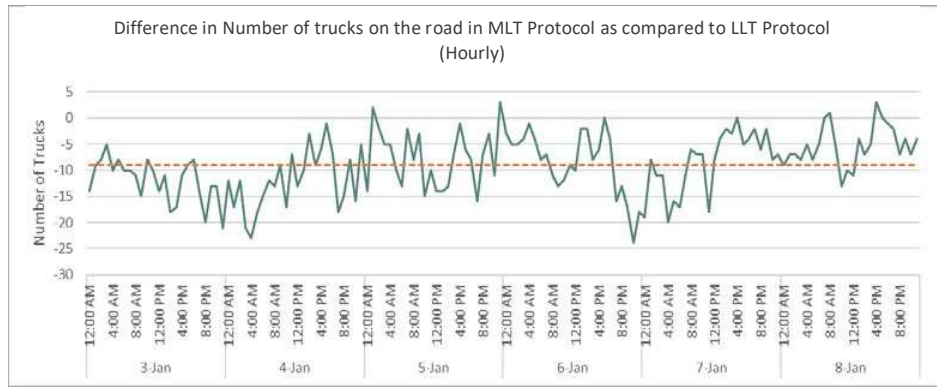

*Figure 3: Difference in the active hourly number of trucks between the two protocols*

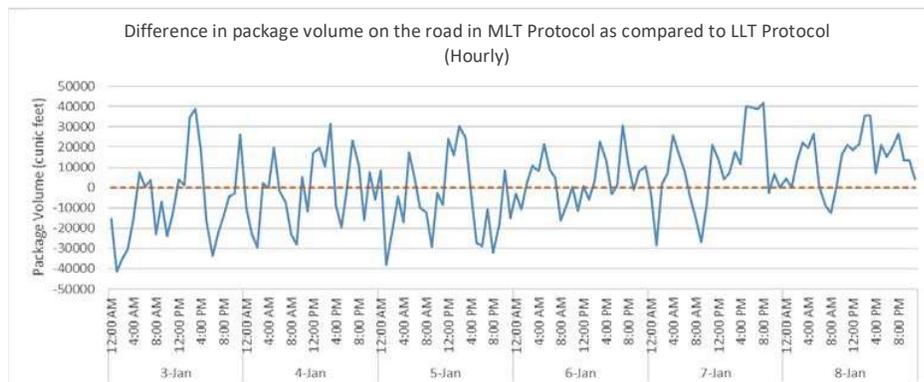

*Figure 4: Difference in the active hourly package volume between the two protocols*

Analyzing the two outputs further, we gain insights on the difference in the number of trucks, total trailer capacity and package volume on the road in the two algorithms as shown in Figure 3 and Figure 4. The total package volume in the MLT protocol is lower than LLT protocol before hiking up notably higher as the former one makes the packages wait at the hubs before they are sent to the next hub. We also observe that in this particular instance, although the number of vehicles in both protocols vary, on average we eliminate 9 trucks at any point in time using MLT protocol as compared to the LLT protocol.

**5.3. Expected Frequency of On-Demand Vehicles**
The simulation model stores and outputs the on-demand vehicle movements, identifying the origin, destination, the departure time and arrival time of each of the used vehicle along with the vehicle capacity and the capacity of the packages that travelled on that vehicle. These outputs can simply be translated into the expected frequency and sizes of vehicles required for each of the segments. In the current simulation results, we notice an overall trend for each of the segments. High flow lanes require a higher number of larger trucks while the low flow lanes have less frequent usage and mostly use smaller trucks. As there is a demand generator built within the simulation model and the package routes are an input to the simulation, these can be altered to obtain different expected frequencies in order to negotiate the on-demand contracts in a better way. The simulation model results for expected frequency will change if the demand distribution input is changed, therefore the simulation model provides an opportunity to perform scenario analysis with different distributions of demand, not just in terms on the quantity but also in terms of the origin and destination of the packages.

## 6. Current Contracts vs. Proposed On-Demand Models

In this section, we compare the proposed on-demand transportation model to the current model used by the national intercity express parcel logistic network. Using the historical data of the vehicle movements for the same demand distribution, we observe that the company is using dedicated fleet as well as outsourced contracts and occasionally obtains spot rates. The dedicated fleet usually travels long haul while the outsourced contracts may be longer segments but these can be one-way trips, and not necessarily a return journey. The spot rates obtained are for one-way trips as well. We compare the active trailer volume on the road in the current contracts and the proposed contracts in 5. In 6, we compare the number of trucks on the road at any hour during the considered time horizon in both models.

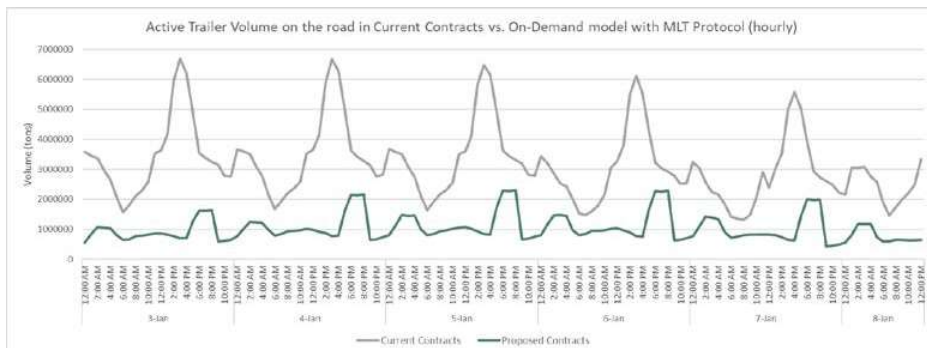

*Figure 5 : Hourly Distribution of Active Trailer Volume*

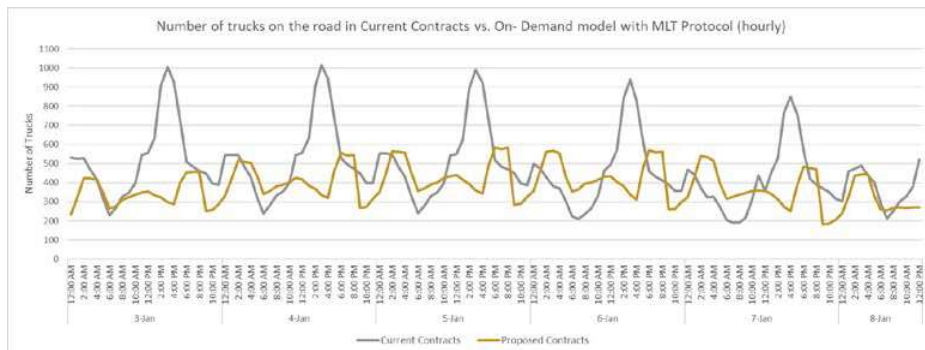

*Figure 6: Hourly Distribution of Number of trucks on the road*

We observe that the total number of vehicle movements has reduced in the on-demand contracts and the total trailer volume has also reduced. As the current contracts are made well in advance without information on the future behaviour of the system, there are many instances where a vehicle may only be 10% filled but still has to depart the hub to be able to reach the next hub within the scheduled time. The flexibility introduced by the on-demand contracts take advantage of the real-time data of packages such as their remaining dwell time to make an informed decision about the capacity and scheduling of the vehicle. With such a model, we are able to reduce the carbon footprint as the number of vehicle movements and total travelled miles required are reduced and are able to achieve a space-efficient solution in terms of higher vehicle fill rates while providing better work environment for drivers as well where they may not be too far from their home but rather shuttle back and forth between two locations.

## 7. Conclusion

The Physical Internet is an open hyperconnected logistics system that uses a set of protocols and interfaces to move physical goods contained within standard modular containers through an open mobility web. In this paper, we have introduced algorithmic protocols to create consolidated sets of modular containers and request the on-demand means to move such sets in the open mobility web, where each container encapsulates packages whose delivery is time sensitive.

We introduced two protocols that enable the scheduling of transportation requests to the contracted service providers, given the assigned dwell time and service promise for packages embedded in containers. We have provided preliminary performance assessment of the protocols using a simulation experiment on actual data of a large parcel delivery service provider. The parameters of the protocol algorithms can be altered according to the network topology as well as the flow and movement of the physical goods. As an example, some low-flow lanes can have a higher Maximum Latency Time ratio, i.e., the containers may wait longer for consolidation versus high flow lanes where the ratio can be set lower.

The Physical Internet requires an extensive set of robust protocols that can be applied to any setting. This paper has focused on one of the facets of such protocols. We believe that our introduction of the proposed protocols may open research avenues, including dynamically optimizing the algorithmic parameters , defining more setting-based parameters to improve the algorithm performance, and performing more extensive simulation-based assessments with scenarios varying notably in terms of demand mix, patterns, and uncertainty; vehicle mix; package size mix; transit time stochasticity; hub availability, capability, and capacity; and tightness of promised delivery times.